\documentclass[11pt]{elsart} 
\usepackage[latin1]{inputenc}
\usepackage[english,swedish,german]{babel}

\usepackage[all,light,draft]{draftcopy}
\pdfoutput=1
\usepackage{graphicx}
\usepackage{psfrag}
\usepackage{amsmath}
\usepackage{tabularx}
\usepackage{verbatim}
\usepackage{float}
\usepackage{amssymb}
\usepackage{latexsym}
\usepackage{pifont}
\usepackage{subfigure}

\newcommand{\vect}{{\operatorname {vec}~}}
\newcommand{\sign}{{\operatorname {sign}}}

\newcommand{\conjj}[1]{{\bar{#1}}}
\newcommand{\re}[1]{{{\mathrm{Re}}~{#1}}}

\newcommand{\acos}{{{\mathrm{acos}}}}
\newcommand{\atan}{{{\mathrm{atan}}}}

\newcommand{\figref}[1]{\ref{#1}}

\newcommand{\RR}{\mathbb{R}}

\newcommand{\CC}{\mathbb{C}}
\newcommand{\LL}{\mathbb{L}}
\newcommand{\MM}{\mathbb{M}}
\newcommand{\ZZ}{\mathbb{Z}}

\newcommand{\NN}{\mathbb{N}}

\newcommand{\kron}{\otimes}

\newcommand{\Arg}[1]{{\mathrm{Arg}~{#1} }}
\newcommand{\LLp}{\LL}

\newcommand{\hp}{\bar{h}}
\newcommand{\fig}[1]{Fig. \figref{#1}}
   {\noindent\begin{description}\item[\fbox{Problem #1}]}
   {\end{description}}

\newenvironment{proof}
   {\begin{pf*}{Proof.}}
   {\hfill $\Box$\end{pf*}}

\newenvironment{example}[1]
   {\begin{exmp}[#1]}
   {\end{exmp}}

\begin{document}

\begin{frontmatter}
\selectlanguage{english}



\title{Critical Delays and 
 Polynomial Eigenvalue Problems}

\author[1]{Elias Jarlebring}

\address[1]{Technische Universität Braunschweig, Institut Computational Mathematics, 38023
 Braunschweig, Germany}
\ead{e.jarlebring@tu-bs.de}


\begin{abstract}
In this work we present a new method to compute the delays of {\em delay 
 differential equations} (DDEs), such that the DDE has a purely imaginary
 eigenvalue. For delay differential equations with multiple delays, the
 {\em critical curves} or {\em critical surfaces} in delay space (that
 is, the set of delays where the DDE has a 
 purely imaginary  eigenvalue) are parameterized. We show how the method
 is related to other  works in the field by treating the case where the
 delays are integer multiples of some delay value, i.e., commensurate
 delays.  

The parametrization is done by solving a {\em quadratic eigenvalue problem} which
 is constructed from the vectorization of a matrix equation and hence 
 typically of large size.  For commensurate delay differential
 equations, the corresponding equation is a polynomial eigenvalue
 problem. As a special case of the proposed method, we find a closed
 form for a parameterization of the  critical surface for the scalar case.  

We provide several examples with visualizations where the computation is
 done with some exploitation of the
 structure of eigenvalue problems.
\end{abstract}

\end{frontmatter}

\section{Introduction}
\selectlanguage{english}
Some phenomena in  engineering, physics and biology can be accurately
described with {\em delay-differential equations} (DDEs), sometimes
referred to as {\em time-delay systems} (TDS). The delay
in these models typically stem from modelling
of phenomena like  transmission, transportation and inertia.  For
instance, in
engineering it is often desirable to design a controller for a system
where physical 
limitations makes the current state of the system
unavailable for measurement. Instead, the feedback in the controller is
done with an old state of the system. This is a typical case where a
delay-differential equation may be used to model the behavior. See
\cite{Niculescu:2001:DELAY} for more examples.  In this paper we
consider linear, time-invariant, $n$-dimensional delay-differential
equations with 
$m$ delays,
\begin{gather}
  \label{intro:tds} 
 \left\{
  \begin{split}
  \dot{x}(t)&=A_0x(t)+\sum_{k=1}^mA_kx(t-h_k),\, t>0\\
  x(t)&=\varphi(t),\, t\in[-h_m,0]
  \end{split}
\right.
\end{gather}
with $x:[-h_m,\infty)\mapsto\RR^{n}$ and 
$h_k\in\RR_+$, $A_k\in\RR^{n\times n}$ for $k=0,\ldots,m$.

For most applications it is desirable that the solution is stable
independent of the initial condition, i.e., that the solution goes
asymptotically to zero when time goes to infinity. This and many other
property can be  determined from
the solutions of the characteristic equation. The {\em characteristic equation} of \eqref{intro:tds} is
\begin{equation}
  \MM(s)v:=\left(-sI_n+A_0+\sum_{k=1}^mA_ke^{-h_ks}\right)v=0, \, \|{v}\|=1,
	\label{intro:char}
\end{equation}
where $v\in\CC^n$ is called {\em eigenvector} and $s\in\CC$ an {\em
 eigenvalue}. The set of all eigenvalues is called the
 {\em spectrum}. 

Similar to ordinary differential equations (ODEs), a DDE is exponentially stable
if and only if all eigenvalues lie in the open left complex half-plane,
because the only clustering point of the real part is minus infinity. An essential difference is that,
unlike ODEs, the spectrum of DDEs contains a
countably infinite number of eigenvalues. Fortunately, it can be
proven (c.f. \cite{Hale:1993:FDE}) that there are only a finite number
of eigenvalues to the right of any vertical line in the complex plane,
e.g., there are only finite number of unstable eigenvalues.

It is of particular interest to determine for what choices of the
delays $h_1, \ldots, h_m$ the DDE \eqref{intro:tds} is stable. This set
of $h_1,\ldots,h_m$ is referred to as the {\em stability region} in
{\em delay-parameter space}. 
From continuity it is clear that, at the boundary of any stability region in
the delay-parameter space, the DDE has at least
one purely imaginary eigenvalue. These critical delay values (single
delay), critical curves (two delays) or critical surfaces (more delays)
are hence important for a complete stability analysis\footnote{The
critical curves (surfaces) are called {\em offspring curves}
and {\em kernel curves} in \cite{Sipahi:2005:COMPLETE}, {\em Hopf bifurcation
curves (surfaces)} in \cite{Hale:1991:GEOMETRIC} and {\em crossing delays
(curves)} in \cite{Malakhovski:2006:2NDORDER} and
\cite{Gu:2005:CROSSING}.}. In this  
work we focus on the computation of critical curves
and surfaces.

In Section~\ref{sect:freedelays} we present a new method to parameterize
such curves and surfaces. An important step of the method is to compute
the eigenvalues of a large {\em quadratic eigenvalue problem},
c.f. \cite{Tisseur:2001:QUADRATIC}, such that the eigenvalue is on the
unit circle. The
matrices of the eigenvalue problem are of size $n^2\times n^2$ and hence
large even for  DDEs of moderate size. The size of the matrices will cause
a computational bottleneck and we therefore make some initial notes on
how an implementation can make use of the structure of the problem of
the corresponding linearized eigenvalue problem. In
Section~\ref{sect:commensuratedelays} we provide a new interpretation of
the results by Chen, Gu and Nett \cite{Chen:1995:DELAYMARGINS}, where 
the {\em commensurate case} is treated. The critical
delays for a system with {\em commensurate delays} can be computed from
the eigenvalues of a polynomial eigenvalue problem of degree $2m$.


The large number of approaches in the literature yielding
delay-dependent stability conditions 
can be classified into two main branches. Most results are either based
on the construction of a Lyapunov-Krasovskii functional, or (as here) based on
a discussion of the imaginary eigenvalues.

Both types of approaches have advantages and disadvantages. For instance, the
methods in \cite{Fridman:2006:NEWCOMPLETE},
\cite{Niculescu:2001:DELAYTRANS}
and the methods summarized in \cite{Gu:2003:STABILITY} are examples
where a Lyapunov-Krasovskii functional is used to construct sufficient
but conservative delay-dependent stability conditions formulated as {\em
linear matrix inequalities} (LMIs). An advantage of these types of
approaches is that the conditions can be elegantly formulated with LMIs
which allow the automatic application in engineering software. A
disadvantage is that the results can be conservative and may
impose unnecessarily strong constraints. Moreover, the treatment of LMIs
of large dimension may be computationally cumbersome. 

The approaches based on imaginary eigenvalues typically yield exact results
but can not be compactly formulated. Early works (and some more recent
works) with a discussion of the imaginary eigenvalues were limited to special
cases. For instance, the trigonometric conditions by Nussbaum \cite[Chapter
III]{Nussbaum:1978:LAGS} are for scalar problems with only two delays,
the conditions by Cooke and Grossman \cite{Cooke:1982:DISCRETE} for
second order DDEs with a single delay, the geometric characterization by
Hale and Huang in \cite{Hale:1993:GLOBAL} for 
first order (scalar) systems with two delays, the discussion using the
quasi-polynomial by Beretta and Kuang in \cite{Beretta:2002:GEOMETRIC}
for single-delay DDEs where the coefficients are dependent on the
delay, the analysis of two-delay DDEs of a special form in
\cite{Gu:2005:CROSSING} (further discussed in Example~\ref{exmp:gnj05}),
the analysis of third order DDEs with two delays by 
Sipahi, et al in \cite{Sipahi:2005:COMPLETE} based on the earlier methods by
Thowsen \cite{Thowsen:1981:ROUTH}, Rekasius
\cite{Rekasius:1980:REKASIUS} and Hertz \cite{Hertz:1984:SIMPLIFIED}, the analysis based on a Pontryagin's
results for quasipolynomials for second order (scalar) DDEs by Cahlon et al in
\cite{Cahlon:2006:STABILITY}, the explicit trigonometric analysis of
scalar two-delay DDEs by Bélair et al in
\cite{Belair:1994:STABILITY}. From this large (but incomplete) list of
works, we conclude that an analysis based on imaginary eigenvalues
is indeed accepted as a natural way to construct  stability conditions.

Recent focus in this field is towards neutral DDEs, distributed DDEs and conditions
which can be efficiently and accurately computed on a computer. The
importance of the possibility to apply the conditions using a computer
should not be underestimated. In fact, the conditions from
most of the methods mentioned above are implicitly or explicitly
formulated using the 
coefficients in the characteristic polynomial with exponential terms
(sometimes referred to as a quasipolynomial). It is well known
from numerical analysis of eigenvalue problem, that this (nominal)
representation of the characteristic equation is not 
(numerically) stable, i.e., the eigenvalues can be (and are often) very
sensitive with respect to 
perturbations (e.g. rounding errors) in the coefficients. This is
important, since it is not possible to represent a polynomial on a
computer using its coefficients without introducing small errors. For this
reason, the methods mentioned above applied to systems of order (say) 4 or more
typically do not give numerical results which can be used in practice. This motivates
the recent works on imaginary eigenvalues using a formulation of the
matrices in the characteristic equation which is generally believed to
scale better with the dimension of the DDE. The
following works represents the class of conditions formulated as 
generalized  eigenvalue problems, Chen, et
al \cite{Chen:1995:DELAYMARGINS} for commensurate 
systems and derivative works \cite{Fu:2006:NEUTRAL,Niculescu:2006:DDAE},
Louisell \cite{Louisell:2001:IMAG} for single delay neutral
systems. In the conditions of these matrix-pencil methods, the
eigenvalues of a large generalized eigenvalue problem must be
computed. Similar to the method discussed here, the matrices are
constructed using the Kronecker product.

We also mention the related work by Ergenc, et al
\cite{Ergenc:2007:EXTENDED} which is partially based on reason with
matrices, but also uses the coefficients in the
quasipolynomial. The method presented here has similarities with this work.


Finally, we note that there are methods to numerically determine parts
of the spectra, e.g. using a {\em linear multistep discretization} of
the solution operator \cite{Verheyden:2007:EFFICIENT} or discretization of the
infinitesimal generator of the semigroup corresponding to the DDE
\cite{Breda:2004:CRDDE}. These methods can be applied for a grid of
points in delay-space. If the grid is fine enough, a  numerical
stability condition is computed.

Apart from the standard works in the field
\cite{Gu:2003:STABILITY,Niculescu:2001:DELAY}, we refer the reader to the
references in the introduction of 
\cite{Knopse:2006:INTERVAL} for a well balanced  overview of modern
delay-dependent stability conditions.



\section{Results}
In order to clearly state the results we start by presenting a
motivating example in Section~\ref{sect:motivation}. The example describes the problem
as well as the idea behind the method. 

The main result consists of the method to parameterize the critical curves
(surfaces) in Section~\ref{sect:freedelays} and some further
interpretation and relation to other publications for the case that the
delays are 
commensurate in Section~\ref{sect:commensuratedelays}. The
computationally dominating part of both of the methods is the
determination of unit eigenvalues of {\em quadratic} or {\em polynomial
eigenvalue problems} of large dimension.

\subsection{Motivation} \label{sect:motivation}
We now describe the general idea behind the method presented in Section~\ref{sect:freedelays}. We do
this by first considering the scalar two-delay case and
describe how the derivation must be changed to hold for multiple
dimensions and multiple delays.

We wish to find purely imaginary eigenvalues, say $s=i\omega$ with
corresponding eigenvector $v$, i.e., $\MM(i\omega)v=0$. In order to
clearly motivate the results we discuss the two delay case and
generalize and formalize the discussion in the following subsections.
The characteristic equation for a DDE with two delays is
\[
   \MM(i\omega)v=\left(-i\omega I+A_0+A_1e^{-ih_1\omega}+A_2e^{-ih_2\omega}  \right)v=0.
\]
Except for some special cases, the points (in delay space) of interest, i.e., points which
causes the DDE to have imaginary eigenvalues, correspond to curves. We wish to
parameterize these {\em critical curves}, and pick
$\varphi:=h_1\omega$ 
as a free parameter. By treating $\varphi$ as a known value we form a
condition on $h_2$ and $\omega$ (such that we can also compute $h_1$).
We will see later that this choice of free parameter is valid for
most cases. The characteristic equation is $2\pi$-periodic in 
$\varphi$, and it will be enough to let
$\varphi$ run the whole span $[-\pi,\pi]$, i.e. for each choice of $\varphi$
in this interval we will find some points  on the critical curves and if
we let $\varphi$ run the whole interval we will find all critical points. The characteristic equation is
such that we wish to find all $\omega\in\RR$ and
$z\in\partial D$ such that
\begin{equation}
   \MM(i\omega)v=\left(-i\omega I+A_0+A_1e^{-i\varphi}+A_2z
					  \right)v=0,\; v\neq 0\label{eq:motphiz}
\end{equation}
where $\partial D$ denotes the unit circle.
If we first consider the scalar case, i.e.,
$A_0=a_0,A_1=a_1, A_2=a_2\in\RR$ the equation corresponds to two
scalar conditions (say real and imaginary parts) and we can eliminate either $\omega$
or $z$. In the approach presented here, we eliminate $\omega$ by forming the
sum of \eqref{eq:motphiz} and its complex conjugate, i.e., 
\[
 0=2a_0+a_1e^{-i\varphi}+a_2\bar{z}+a_1e^{i\varphi}+a_2z=
a_2\bar{z}+2a_0+2a_1\cos(\varphi)+a_2z.
\]
Multiplying with $z$ yields the quadratic equation,
\[
 a_2z^2+2z(a_0+a_1\cos(\varphi))+a_2=0,
\]
since $z\bar{z}=1$.
It has the two solutions
\[
 z=
\frac{ -(a_0+a_1\cos(\varphi))\pm i\sqrt{a_2^2-(a_0+a_1\cos(\varphi))^2}}{a_2},
\]
assuming $a_2\neq 0$.
We can now compute $\omega$ by inserting $z$ into \eqref{eq:motphiz} and
rearrange the terms, i.e.,
\[
 i\omega=a_0+a_1e^{-i\varphi}+a_2z=i\left(-a_1\sin(\varphi)\pm\sqrt{a_2^2-(a_0+a_1\cos(\varphi))^2}\right).
\]
Since $z=e^{-ih_2\omega}$ and $\varphi=h_1\omega$, a parametrization of
the critical curves is
\begin{equation}\label{eq:motivatingparametrization}
 \hp(\varphi)=
\begin{pmatrix}
h_1\\h_2
\end{pmatrix}=
\begin{pmatrix}
\frac{\varphi+2p\pi}
{-a_1\sin(\varphi)\pm\sqrt{a_2^2-(a_0+a_1\cos(\varphi))^2}}
\\
\frac{-\Arg z+2q\pi}
{-a_1\sin(\varphi)\pm\sqrt{a_2^2-(a_0+a_1\cos(\varphi))^2}}
\end{pmatrix},
\end{equation}
where
$\Arg z=\pm\sign(a_2)\acos\left(-\frac{a_0+a_1\cos(\varphi)}{a_2}\right)$
for any $p,q\in\ZZ$ and $\varphi\in[-\pi,\pi]$. Note that the signs must
be matched, i.e., for each choice of the free parameter $\varphi$ and
branches $p$ and $q$, the parameterization has two critical delays.

There are other approaches based on other parameterizations. For
instance, the analysis of Gu et
al \cite{Gu:2005:CROSSING} can be seen as a parameterization with
$\omega\in\RR$ as free parameter. They 
reach the result that if the characteristic equation can be rewritten to the form
\[
 1+a_1(s)e^{-h_1s}+a_2(s)e^{-h_2s}=0,
\]
where $a_1$ and $a_2$ are rational functions, a parameterization of
the critical curves are given by
\begin{eqnarray*}
h_1&=&\frac{\Arg (a_1(i\omega)) +(2p-1)\pi\pm \theta_1}{\omega},\;\; \theta_1=\acos\left(\frac{1+|a_1(i\omega)|^2-|a_2(i\omega)|^2}{2|a_1(i\omega)|}\right),\\
h_2&=&\frac{\Arg (a_2(i\omega)) +(2q-1)\pi\mp\theta_2}{\omega},\;\;\theta_2=\acos\left(\frac{1+|a_2(i\omega)|^2-|a_1(i\omega)|^2}{2|a_2(i\omega)|}\right).
\end{eqnarray*}

We also note that the different types of possible critical curves and
 other properties are
classified in \cite{Hale:1993:GLOBAL}.

\begin{example}{Different parameterizations}
\label{exmp:gnj05}
 In order to show the difference between
 \eqref{eq:motivatingparametrization} and the parameterization in \cite{Gu:2005:CROSSING},
 we construct the parameterizations of the
 critical curves for the case that  $a_0=a_1=-1$,
 $a_2=-\frac12$. According to \eqref{eq:motivatingparametrization}, the critical curves are given by,
\[
 h(\varphi)=
\begin{pmatrix}
\frac{\varphi+2p\pi}
{\sin(\varphi)\pm\frac12\sqrt{5+8\cos(\varphi)+4\cos^2(\varphi)}}
\\
\frac{\pm\acos\left(-2-2\cos(\varphi)\right)+2q\pi}
{\sin(\varphi)\pm\frac12\sqrt{5+8\cos(\varphi)+4\cos^2(\varphi)}}
\end{pmatrix}.
\]
In this example we can find a slightly more elegant
 parametrization by letting $x=\cos(\varphi)$, i.e.,
\[
 h(x)=
\begin{pmatrix}
\frac{\acos(x)+2p\pi}
{-\sign(x)\sqrt{1-x^2}\pm\sqrt{0.25-(1+x)^2}}\\
\frac{\pm\acos(-2-2x)+2q\pi}
{-\sign(x)\sqrt{1-x^2}\pm\sqrt{0.25-(1+x)^2}}
\end{pmatrix}, x\in[-1,1].
\]
The critical curves are shown in Figure~\ref{fig:gnj05}. The minimum of
 the 2-norm, i.e. the 2-norm stability delay-radius, is $\|h(x)\|_2\approx 2.896$ and taken at
 $x\approx -0.7012$ where $\omega\approx1.1139$, $h_1\approx 2.1078$ and
 $h_2\approx 1.9853$.
\begin{figure}[h]
  \centering
  {\scalebox{0.6}{\includegraphics{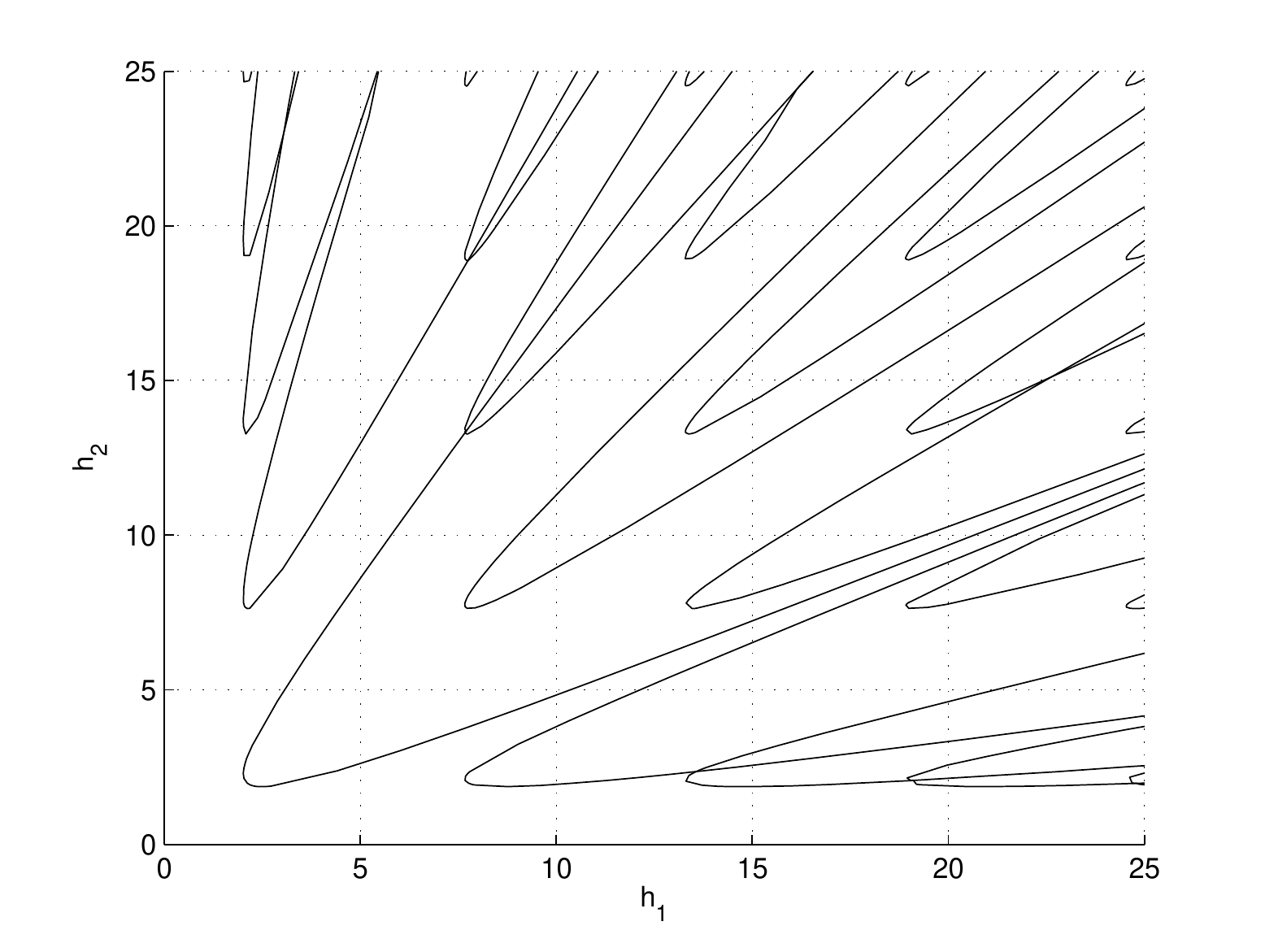}}}
  \caption{Critical curves for Example~\ref{exmp:gnj05}}
 \label{fig:gnj05} 
\end{figure}

In the context of \cite{Gu:2005:CROSSING},
 $a_1(i\omega)=\frac1{1+i\omega}, a_2(i\omega)=\frac1{2(1+i\omega)}$ and the
 parameterization is given by
\begin{eqnarray*}
   h_1&=&\frac{-\Arg (1+i\omega)+ (2p-1)\pi \pm\theta_1(\omega)}{\omega} =\\
&=&\frac{1}{\omega}\left(-\atan(\omega)\pm\acos\left(\frac12\sqrt{1+\omega^2}+\frac3{8\sqrt{1+\omega^2}}\right)+(2p-1)\pi\right)\\
   h_2&=&\frac{-\Arg (1+i\omega)+ (2q-1)\pi
	 \mp\theta_2(\omega)}{\omega}\\
&=&\frac{1}{\omega}\left(-\atan(\omega)\mp\acos\left(\sqrt{1+\omega^2}-\frac3{4\sqrt{1+\omega^2}} \right)+(2q-1)\pi\right),
\end{eqnarray*}
which represent the same set of curves, but has a surprisingly small  obvious similarity
 with the other parameterization \eqref{eq:motivatingparametrization}. 
\end{example}



We now wish to do an analysis similar to the derivation of \eqref{eq:motivatingparametrization} for the general case where
$A_0,A_1$ and $A_2$ are not scalar. It turns out that a straightforward
generalization does not yield such explicit results. The reason is that
the elimination of $\omega$ must be done in a slightly different
way. The first step in the discussion above involves taking the sum of
\eqref{eq:motphiz} and its complex conjugate. Since $v$ is typically
a complex vector, the sum of \eqref{eq:motphiz} and its  complex conjugate 
\eqref{eq:motphiz} will not eliminate $\omega$. Taking the sum of the
\eqref{eq:motphiz} and its complex conjugate transpose is of course also
not possible as the dimensions do not fit. Instead we take the sum of
$\MM(s)vv^*$ and $vv^*\MM(s)^*$, where $A^*$ denotes the complex
conjugate transpose of $A$. We believe this is a natural way to generalize the
elimination of $\omega$.
Clearly,
\begin{multline*}
  0=\MM(s)vv^*+vv^*\MM(s)^*= \\
  =A_2vv^*z+((A_0+A_1e^{-i\varphi})vv^*+vv^*(A_0^T+A_1^Te^{i\varphi}))+vv^*A_2^T\bar{z},
\end{multline*}
which is a matrix equation. This equation can be rewritten into
a {\em quadratic eigenvalue problem} by vectorization, i.e., stacking the
columns of the matrix on top of each other. Quadratic eigenvalue
problems can be rewritten into a generalized eigenvalue
problem using a so called {\em companion linearization} (c.f
\cite{Mackey:2006:VECT}). 
As we will see later, if $n$ is not large, there are  
numerical methods to find all $z$. We can then compute the corresponding
$\omega$ and the critical delays similar to the way done for the scalar case.

Even though the motivation above is simple, it is not clear
that all steps involved are reversible, i.e., there is no convincing
argument to ensure  that we get
all critical delays. For that reason, we formalize the discussion in the
sections that follow. At the same time
we generalize the method to an arbitrary number of delays. The
generalization is such that if we have $m$ delays we have $m-1$ free
parameters, and need to solve a quadratic eigenvalue problem for each
choice of the parameters. In order to compare
the method with related approaches we also discuss
the  case where the delays are (fixed) integer multiples of the first
delay, i.e., commensurate delays.

\subsection{Main results}\label{sect:mainresults}

Our main goal is to formalize and generalize the discussion in the
previous section, and to see how the constructed method relates to
approaches in the literature. In particular, we wish to provide further understanding to the method in 
\cite{Chen:1995:DELAYMARGINS} where a matrix-pencil method is presented
for the {\em commensurate case}, i.e., the case where the delays are integer
multiples of the first delay. We do this by first, in a general setting,
introducing an operator $\LL$. For this operator  we can formulate an
equivalence theorem with the characteristic equation. 
One reason for formulating such a theorem 
 is that the results in \cite{Chen:1995:DELAYMARGINS} are easily interpreted from the
conditions in the equivalence theorem. Another reason is that from
the conditions we can construct a new method (similar to the motivation)
to compute the
critical delays for the case that the delays are not necessarily commensurate. Hence, both results can be
interpreted in terms of the conditions in the equivalence theorem. 


\begin{thm}\label{thm:equiv}
  Given $s\in \CC$ and $v\in \CC^n$, $v^*v=1$ the following statements are equivalent.
\begin{eqnarray}
		  &&\MM(s)v=0 \label{equiv:1} \\
 && \LLp(vv^*,s)=0\, \wedge\, v^*\MM(s)v=0  \label{equiv:2}
\end{eqnarray}
where
\begin{gather}
\begin{split}
	 \LLp&(X,s):=\left(\MM(s)\right)X+X\left(\MM(s)^*\right)=\\
     &=\sum_{k=0}^{m}\left(A_kXe^{-h_ks}+XA_k^Te^{-h_k\conjj{s}}
 \right)-2X\re{s},
\end{split}
\label{def:lyap}
\end{gather}
and $h_0=0$ for notational convenience.
\end{thm}
\begin{proof}
  The forward implication is trivial from definitions, i.e., \eqref{intro:char}
 and \eqref{def:lyap}. The backward
 implication, i.e., \eqref{equiv:2}$\Rightarrow$\eqref{equiv:1}, is clear from the following equality.
\[
  \MM(s)v =\LLp(vv^*,s)v-vv^*\MM(s)^*v
\]
\end{proof}
The reason why we construct $\LLp$ in this way, is that similar to the
motivation, for the critical case, i.e.,  $s=i\omega\in i\RR$,  $s=i\omega$ only appears in
the exponential terms of the expression.

In the two following sections we rephrase the equivalence theorem such
that we get a method to compute the critical delays. Because of the
equivalence theorem the formulation will not miss any critical delays.

\subsection{Free delays}\label{sect:freedelays}

In the motivation we saw that for two delays  we needed to parameterize
a {\em curve} in delay-space. In general, if we have
$m$ delays we need to parameterize the $(m-1)$ dimensional
(hyper-)surface. We
do that by letting $z=e^{-ih_m\omega}$ and represent the other exponential
terms by free variables, i.e., $\varphi_k=h_k\omega$,
$k=1,\ldots,m-1$. This allows us to rewrite
the main expression of  Theorem~\ref{thm:equiv} to a
quadratic eigenvalue problem which can be solved by transforming it to
an eigenvalue problem using a {\em companion linearization}. 

For the scalar case, the quadratic eigenvalue problem can be solved
explicitly, and we can form an explicit parametrization of the critical
surfaces. We also verify the results by comparing it to classical
results for the single delay, scalar delay-differential equations.

The following theorem follows from the substitution 
$z=e^{-ih_m\omega}$,$\varphi_k=h_k\omega$, $k=1,\ldots,m-1$
and  Theorem~\ref{thm:equiv}.
\begin{thm}
  The critical delays are given by
 $h_m=\frac{-\Arg{z}+2p_m\pi}{\omega}$,
 $h_k=\frac{\varphi_k+2p_k\pi}{\omega}$, $k=1,\ldots,m-1$, where 
$p_k\in \ZZ, k=1\ldots m$, $\varphi_k\in[-\pi,\pi], k=1,\ldots,m-1$,
 $z\in\partial D$, $\omega\in \RR$
 such that 
\begin{equation}
z^2 A_mvv^* +z\left(\sum_{k=0}^{m-1}A_kvv^*e^{-i\varphi_k}+vv^*A_k^Te^{i\varphi_k}\right)
  +vv^*A_m^T=0, \label{eq:critquad}
\end{equation}
\begin{equation}
 \omega=-iv^*\left(A_mz+\sum_{k=0}^{m-1}A_ke^{-i\varphi_k}\right)v,\label{eq:critomega}
\end{equation}
where $\varphi_0=0$ for notational convenience.
\label{thm:quadthm}
\end{thm}



Apart from singular cases (which will be discussed later), this theorem
generates a parameterization of the critical
curves, with $\varphi_1,\ldots,\varphi_{m-1}$ as free variables. It is clear that if we can find solutions $z$, $v$ of
\eqref{eq:critquad}, $\omega$ and critical delays can be computed for a
specific choice of $\varphi_k$, $k=1,\ldots,m-1$. 
Hence, the difficulty of applying this theorem is  solving
\eqref{eq:critquad}. We solve this by vectorizing and rewriting the
quadratic equation to a generalized eigenvalue problem.


The vectorized version of \eqref{eq:critquad}, i.e., the same
equation but the columns of the matrix condition stacked on each other,
is the equation
\begin{equation}
\left(z^2 I\otimes A_m +z\sum_{k=0}^{m-1}L_k(\varphi_k)
  +A_m\otimes I\right)u=0, \label{results:critquadvect}
\end{equation}
where  $u=v\otimes\bar{v}$ is the vectorization of $vv^*$ and $L_k(\varphi_k):=I\kron
 A_ke^{-i\varphi_k}+A_k\kron I e^{i\varphi_k}$, and
$\otimes$ the Kronecker product with the usual meaning. Problems of the type
 $(Mz^2+Gz+K)v=0$, e.g., \eqref{results:critquadvect},  are  known in the
 literature as  {\em quadratic eigenvalue problems}
 (c.f. \cite{Tisseur:2001:QUADRATIC}). They are typically solved by
 transformation into first order form, normally using a so called
 {\em companion linearization}, e.g., here
\begin{equation}
  \begin{pmatrix}
0& I\\
I\kron A_m&    \sum_{k=0}^{m-1}L_k(\varphi_k) \\
\end{pmatrix}\begin{pmatrix}u\\z u\end{pmatrix}=
  z\begin{pmatrix}
I&0\\
0 &-A_m\kron I\\
\end{pmatrix}\begin{pmatrix}u\\z u\end{pmatrix}.
\label{eq:companionform2}
\end{equation}

We can hence find the solutions of \eqref{eq:critquad} by finding
the eigenvalues corresponding to \eqref{eq:companionform2}, to which
standard methods for generalized eigenvalue problem can be applied. There are
also other numerical methods to solve quadratic eigenvalue problems
without linearizing it,
e.g., {\em Jacobi-Davidson} \cite{Sleijpen:1996:POLYNOMIAL} and {\em
second order Arnoldi} \cite{Bai:2005:SOAR}. In this paper we focus on
the companion linearization \eqref{eq:companionform2}.

With these methods it possible (unless $n$ is large) to find the
eigenpairs of \eqref{eq:companionform2} to sufficient accuracy on a
computer. Because it is not possible (in general) to solve eigenvalue
problems without errors, we will not get exactly a vector corresponding to a
rank-one matrix. However, if the eigenvalue is isolated we can find an
approximation close to the rank-one matrix. We pick the
approximation by choosing the rank-one matrix corresponding to the
singular value which is not close to zero. All but one singular value
should be close 
to zero, because the numerical rank is one.  The accuracy of this choice
of the approximation can be characterized by the following simple result.
\begin{lem}\label{lem:quadaccuracy}
Let $(z,x)$ be an eigenpair of the quadratic eigenvalue problem
 $(z^2M+zG+K)x=0$ and $x=v\otimes \bar{v}$, $\|x\|=1$, $|z|=1$. Say
 $\tilde{x}=x+y=u\otimes 
 \bar{u}+q$ is an approximation of $x$, then the sine of the angle between
 the approximation $u$ and the vector $v$ is bounded by
\[
 |\sin(u,v)|\le \sqrt{\|y\|\left(\|M+I\|+\|G\|+\|K\|\right)+\|q\|}.
\]
\end{lem}
\begin{proof}
Multiplying the characteristic equation from the left with $\tilde{x}^*$
 and adding $z^2\tilde{x}^*x$ to both sides, yields
\[
 \tilde{x}^*xz^2=\tilde{x}^*(z^2(M+I)+zG+K)x.
\]
We exploit that
 $\tilde{x}^*x=(u\otimes\bar{u}+q)^*(v\otimes\bar{v})=u^*v\bar{u}^*\bar{v}+q^*x$
 from which we deduce that
\begin{multline*}
  z^2|u^*v|^2=z^2(u^*v\bar{u}^*\bar{v})=\tilde{x}^*(z^2(M+I)+zG+K)x-z^2q^*x=
\\
=x^*(z^2(M+I)+zG+K)x+y^*(z^2(M+I)+zG+K)x-z^2q^*x=\\
=z^2+y^*(z^2(M+I)+zG+K)x-z^2q^*x,
\end{multline*}
i.e.,
\[
  |u^*v|^2=1+y^*((M+I)+z^{-1}G+Kz^{-2})x-q^*x
\]
and
\[
 |u^*v|^2\ge 1-\|y\|\left(\|M+I\|+\|G\|+\|K\|\right)-\|q\|,
\]
where we used that $|z|=1$. Finally, the result follows from the fact
 that $|\sin(u,v)|=\sqrt{1-|u^*v|^2}$.
\end{proof}

In rough terms, the result is that if the error of the
eigenvector, denoted by $y$, is small and the distance to the Hermitian
rank-one matrix $q$ is small, then the angle between $v$ and its'
approximation $u$ is small.  The lemma also indicates the (unfortunate)
fact that if the error of the eigenvector is $\epsilon_{mach}$ the accuracy
of the approximation $u$ is essentially $\sqrt{\epsilon_{mach}}$.

If we have multiple eigenvalues, the
invariant subspace must be sought for rank-one matrices. This
typically unusual case will not be treated in detail. 

Moreover, we are
only interested in eigenvalues $z$ on the unit circle. Because of
rounding errors, we classify an eigenvalue as a unit-eigenvalue if its
modulus differs from $1$ less than some tolerance, say $\sqrt{\epsilon_{mach}}$.

It is also worth noting that \eqref{results:critquadvect} is a quadratic
eigenvalue problem and can have (so called) {\em infinite
eigenvalues}. For instance, if $A_m=0$ the equation is independent of 
$z\in\CC$, i.e., if the linear term is zero it is fulfilled for
any $z\in\partial D$. The interpretation of 
the theorem is then that the criticality is independent of
$h_m$, and one of the free parameters is no longer free.
Since this can only occur if $A_m$ is singular, the problems of infinite
eigenvalues can be avoided by reordering of the matrices (assuming at
least one matrix is nonsingular). We characterize the phenomena of infinite
eigenvalues with the following example.
\begin{example}{Singular case}\label{exmp:singular}
Consider the DDE with
\[
  A_0=\begin{pmatrix}0&0\\0&0
\end{pmatrix}
,\,A_1=\begin{pmatrix}2&\varepsilon\\3&1\end{pmatrix},
A_2=\begin{pmatrix}1&0\\0&0
\end{pmatrix}.
\]
With this example we wish to show how an infinite eigenvalue of
 \eqref{results:critquadvect} can be interpreted, and how this appears
 in 
 Theorem~\ref{thm:quadthm}. For this
 particular example it is actually possible to circumvent the problem
 with the infinite eigenvalue simply by switching $A_1$ and $A_2$. Here,
 we
 will not approach the problem that way, because that type of
 reordering of the matrices is not always possible. 

 If \eqref{results:critquadvect} has an infinite eigenvalue, it must
 have a corresponding eigenvector such that the quadratic
 term disappears, i.e., here  $u=(0,1)^T\kron (0,1)^T=(0,0,0,1)^T$ is
 the only eigenvector of the correct form. For this $v$, the
 quadratic and the constant term in \eqref{eq:critquad} are zero, and
 \eqref{eq:critomega} is independent of $z$. Since $z\neq 0$, $v$ is a
 valid solution of \eqref{eq:critquad} only if 
\begin{equation}
 A_1vv^*e^{-i\varphi_1}  +vv^*A_1^Te^{i\varphi_1}=0.\label{eq:singularexmp}
\end{equation} This means that
 the only choice of  $\varphi_1$ which generates a critical delay
 (corresponding to this $v$) is when \eqref{eq:singularexmp} is fulfilled.
Here \eqref{eq:singularexmp} and hence \eqref{eq:critquad} is
 fulfilled only if $\varepsilon=0$ and for $\varphi_1=-\frac{\pi}2+2p\pi$ with the
 corresponding critical
 curve $h_1=-\frac{\pi}{2}+2p\pi $ (for any $h_2$). 
The critical curves are plotted for three choices of $\varepsilon$ in
 Figure~\ref{fig:singular}. From the figure it is easy to identify that
 for $\varepsilon=0$ there are a vertical critical curves, i.e., curves which
 are independent of $h_2$ corresponding to the contribution of the
 infinite eigenvalue. If $\varepsilon>0$ there are no vertical lines, 
 no infinite eigenvalues and the critical curves are characterized by the
 finite eigenvalues.
\begin{figure}[h]
  \centering
  {\scalebox{0.6}{\includegraphics{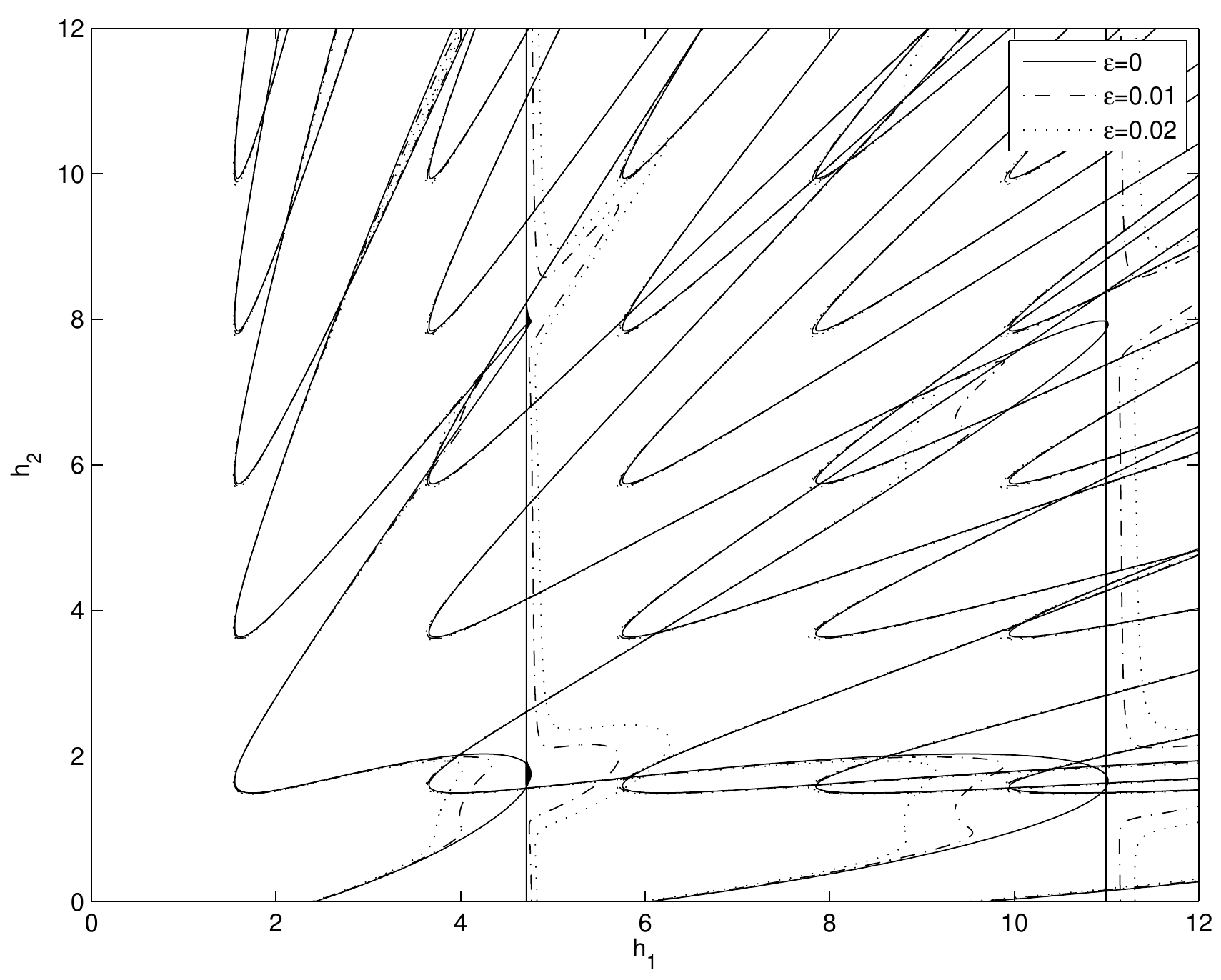}}}
  \caption{Critical curves for Example~\ref{exmp:singular}}
 \label{fig:singular} 
\end{figure}

\end{example}


For the  scalar  case, the  quadratic eigenvalue problem
reduces to a quadratic equation. The
theorem can then be simplified to

\begin{cor}
  For the case that $A_k=a_k\in\RR$, $k=0,\ldots,m$, the critical delays are given by
\begin{eqnarray*}
h_k&=&\frac{\varphi_k+2p_k\pi}{\omega},\,k=1,\ldots,m-1\\
h_m&=&\frac{\mp\sign(a_m)\acos\left(\frac{\sum_{k=0}^ma_k\cos(\varphi_k)
									  }{a_m}\right)+2p_m\pi}{\omega},\\
\omega&=&-\sum_{k=0}^ma_k\sin(\varphi_k)\pm\sqrt{a_m^2-\left(\sum_{k=0}^ma_k\cos
																	(\varphi_k)\right)^2}
\end{eqnarray*}
 where 
$p_k\in \ZZ, k=1\ldots m$, $\varphi_k\in[-\pi,\pi], k=1,\ldots,m-1$.
\label{thm:secondary-onedim}
\end{cor}

It is worth noting that the formulas in the literature for the  scalar
single delay  case normally assume $a_1\le-|a_0|$, because otherwise the
DDE is (un)stable independent of delay or the corresponding delay-free DDE  is not stable. From Corollary~\ref{thm:secondary-onedim},
the critical delays for the single delay scalar DDE
are given by
\[
 h=\frac{-\sign(a_1)\acos(-\frac
 {a_0}{a_1})+2p\pi}{\sqrt{a_1^2-a_0^2}},
\]
which reduces to classical results if $a_1<0$, c.f. \cite{Cooke:1982:DISCRETE} or  \cite[Section
3.4.1]{Niculescu:2001:DELAY}

\subsection{Commensurate delays}\label{sect:commensuratedelays}
For the case that the delays $h_k$ are integer multiples of the smallest
delay, say $h$, Theorem~\ref{thm:equiv} can be restated such that it is
very similar to the results by Chen, Gu and Nett in
\cite{Chen:1995:DELAYMARGINS}. For $h_k=n_kh$, the problem is to find the critical
delay $h$ such that the choice of the delays in the direction
$\hp=(h_1,h_2,\ldots,h_m)=(n_1h,n_2h,\ldots,n_mh)$, where
$n_1,\ldots,n_m\in\NN$, generates a purely imaginary eigenvalue. In
other words, we fix a ``rational'' direction in delay-space with
$(n_1,\ldots,n_m)$  and search for critical delays along this direction.

The following theorem follows from Theorem~\ref{thm:equiv} by letting $z=e^{-ih\omega}$.
\begin{thm}
   The critical delays $h\in\RR_+$ corresponding to the direction in
 delay space defined by
\[
  \hp=(h_1,\ldots,h_m)=(hn_1,hn_2,\ldots,hn_m),
\]
 are given by 
 \[
    h=\frac{-\Arg z +2p\pi}{\omega},
 \]
for any $p\in\ZZ$ and for $v\in \CC^{n}$, $v^*v=1$, $\omega\in
 \RR$, $z\in\partial D$ fulfilling 
	\begin{equation}
	 \sum_{k=0}^{m}\left(A_kvv^* z^{n_k}+vv^*A_k^T z^{-n_k}
							\right)=0,\, \label{eq:chengu}
	\end{equation}
   and
   \[
      i\omega=v^*\left(\sum_{k=0}^mA_k z^{n_k}\right)v,
	\] 
where $n_0=0$ for notational convenience.
 \label{thm:chengu}
\end{thm}
Analogously to the previous section, if it is possible to find the solutions 
$z$ and $v$ of
\eqref{eq:chengu} the other expressions explicitly yield $\omega$ and $h$. 
Without loss of generality we let $n_m=\max_{k\in[1,\ldots,m]} n_k$. After 
vectorizing the matrix equation \eqref{eq:chengu} we find that
\begin{equation}
   \sum_{k=0}^m
\left(I\otimes A_k z^{n_m+n_k}+A_k\otimes I z^{n_m-n_k}
							\right)u=0\label{eq:chengupoly}
\end{equation}
where $u\in\CC^{n^2}$ is the vectorization of $vv^*$, i.e., $u=v\otimes
\bar{v}$. This equation is of the form
\begin{equation}
  \sum_{k=0}^{N}B_kz^ku=0, \label{eq:polyform}
\end{equation}
which in the literature is known as a {\em polynomial eigenvalue problem}.
Similar to quadratic eigenvalue problems, the most common way to solve
polynomial  
eigenvalue problems  is by companion linearization, which is analyzed and
generalized in \cite{Mackey:2006:VECT} and \cite{Mackey:2006:STRUCTURED}.  For instance, the eigenvalues of
\eqref{eq:polyform} are the eigenvalues corresponding to the generalized
eigenvalue problem
\begin{equation} \label{eq:companionform}
   z\begin{pmatrix} 
	  I&        & &      \\
      &\ddots  & &		\\
      &        &I&		\\
      &        & &B_N\\
	\end{pmatrix}w=    
   \begin{pmatrix}
	  0 &  I   &               \\
	    &  \ddots   &  \ddots          \\
	    &      &  0      & I \\
   -B_0& \cdots & -B_{N-2} &-B_{N-1} \\
	\end{pmatrix}w,
\end{equation}
where $w=(u^T, z u^T, z^2u^T,\cdots z^{N-1}u^T)^T$. By 
selecting $B_k$ according to the coefficients in
\eqref{eq:chengupoly}, we can  compute all solutions 
$z, v$
to \eqref{eq:chengu} by solving the eigenvalue problem 
\eqref{eq:companionform}. Other methods, such as the {\em Jacobi-Davidson} method 
could be applied directly to the polynomial eigenvalue problem without
linearizing it, c.f., \cite{Sleijpen:1996:POLYNOMIAL}.

The Hermitian rank-one matrix can be chosen similar to the choice in the
previous section, i.e., the principal vector in the singular value
decomposition, since the accuracy result in Lemma~\ref{lem:quadaccuracy}
generalizes to the polynomial case.
\begin{lem}\label{lem:polyaccuracy}
Let $(z,x)$ be an eigenpair of the polynomial eigenvalue problem
  \eqref{eq:polyform} and $x=v\otimes \bar{v}$, $\|x\|=1$, $|z|=1$. Say
 $\tilde{x}=x+y=u\otimes 
 \bar{u}+q$ is an approximation of $x$, then the sine of the angle between
 the approximation $u$ and the vector $v$ is bounded by
\[
 |\sin(u,v)|\le
 \sqrt{\|y\|\left(1+\sum_{k=0}^{N}\|B_k\|\right)+\|q\|}. 
\]
\end{lem}
\begin{proof}
The proof is analogous to the proof of Lemma~\ref{lem:quadaccuracy}.
\end{proof}

\begin{rem}
 For the  case $h_k=hk$, the companion form
 \eqref{eq:companionform} is very similar to the 
 eigenvalue problem occurring  in \cite{Chen:1995:DELAYMARGINS}. However, in that
 context it is not recognized that the eigenproblem to be solved is a
 polynomial eigenproblem and that the eigenvector  $u$ is the
 vectorization of an Hermitian rank one 
 matrix. 
\end{rem}
\subsection{Notes on computation} \label{sect:numerics}
We now clarify how we can use Theorem \ref{thm:quadthm} to generate the
critical surfaces, discuss computational bottlenecks and suggest
possible improvements and alternatives.

Theorem \ref{thm:quadthm} is an equivalence theorem between
$\varphi_1,\ldots,\varphi_{m-1}$ and $h_{1},\ldots,h_{m}$. Hence, we can
see it as a parameterization of the critical surface. In other words, if
we let $\varphi_1,\ldots,\varphi_{m-1}$ run the whole interval, we will generate all critical
points. In practice, we typically let the free parameter run over finite
number of grid points with a grid size small enough such that we can
convince ourselves of the continuity of the critical curves.

With this in mind, we state in pseudo-code how to
generate the critical curves for the two-delay system.

\newcommand{\myspace}{\hspace{0.3cm}}
\hspace{-0.3cm}
\begin{minipage}{0.7\textwidth}
\begin{enumerate}
\small
\item[1.] FOR $\varphi=-\pi:\Delta:\pi$
\item[2.] \myspace Find eigenpairs $(z_k,u_k)$ of \eqref{eq:companionform2}
\item[3.] \myspace FOR $k=1:\operatorname{length}(z)$
\item[4.] \myspace\myspace IF $z_k$ is on unit circle
\item[5.] \myspace\myspace\myspace Compute $v_k$ such that $u_k=\vect{v_kv_k^*}$
\item[6.] \myspace\myspace\myspace Compute
		$\omega_k=-iv_k^*\left(A_2z_k+A_0+A_1e^{-i\varphi}\right)v_k$
\item[7.] \myspace\myspace\myspace Accept critical points $(h_1,h_2)$
\begin{eqnarray*}
 h_1&=&\frac{\varphi+2p\pi}{\omega_k},\,p=-p_{max},\ldots,p_{max}\\
 h_2&=&\frac{-\Arg{z_k}+2q\pi}{\omega_k},\,q=-p_{max},\ldots,p_{max}
\end{eqnarray*}
		
\item[8.] \myspace\myspace END
\item[9.] \myspace END
\item[10.] END
\end{enumerate}
\end{minipage}

In step 1, $\Delta$ is the stepsize of the parameter $\varphi$. In
step 7, $p_{max}$ is the number of branches which should be
included in the computation. Step 7 is not computationally
demanding. 
We can select $p_{max}$ so large that the
computation contains all relevant branches, say such that all delays
smaller than some delay tolerance are found. This is possible because
the delays are monotonically increasing or decreasing in the branch parameter. The generalization to more than two delays is
straighforward. It involves a nesting of the outer iteration (step 1)
with for-loops of the new free variables $\varphi_k$ and computing the other delays
in step 7 similar to $h_1$.

For DDEs of large or moderate dimension, the main computational effort in an implementation of the method is the solving of eigenvalue problems. This can
be seen from the following remarks on the computational complexity.

The current general full eigenproblem solvers (for instance the QR-algorithm) require, roughly
speaking, a computational effort proportional to $N^3$  to compute the
eigenvalues of a $N\times N$-matrix to sufficient accuracy. 
For Theorem~\ref{thm:quadthm}, the companion matrices are of dimension $N=2n^2$ and
for Theorem~\ref{thm:chengu}, $N=2n_mn^2$. For both cases the
computational effort is hence proportional to $n^6$, which is a lot 
for moderate or large system-dimensions $n$. For sparse eigensolvers the
complexity is reduced, but the dimension of the system still causes a
computational bottleneck.

\begin{rem}
We note that the very high computational cost for large systems, even
 suggests that instead of parameterizing using ${m-1}$ free variables
 over $\varphi\in[-\pi,\pi]^{m-1}$, where each evaluation costs $n^6$, it may
 be more efficient to search $\varphi\in[-\pi,\pi]^m$ and as an
 additional constraint check if 
 (or how far from) an imaginary eigenvalue is contained in the spectrum of
 $A=A_0+\sum_{k=1}^mA_ke^{-i\varphi_k}$. Computing the spectrum of $A$
 has the computational cost $n^3$, which (for large $n$) compensates for
 the additional degree of freedom which has to be scanned. We note that
 some form of iterative method is likely to be necessary in order
 to impose the additional constraint (the imaginary eigenvalue) in a
 reliable way.

 This remark was brought to our attention by Michiel Hochstenbach.
\end{rem}

In order to make larger problems tractable we make some initial remarks
on how to adapt numerical methods for \eqref{eq:companionform2}.
The eigenvalues we wish to find have the following special properties
we can make use of.
\begin{enumerate}
 \item{The eigenvalues of interest lie on the unit circle.}
 \item{The matrices resulting from the companion form are sparse.}
 \item{The matrices have a very special structure originating from the
		Kronecker construction.}
 \item{Only eigenvectors of the form $u=\vect vv^*$, i.e. a
		vectorization of an hermitian rank one
		matrix, are of interest.}
 \item{The eigenvalues move continuously with respect to
		the parameters.}
\end{enumerate}

In this work we only exploit the first two items. The unit eigenvalue
property is exploited by transforming the unit eigenvalues to the real line,
i.e. we make the Cayley-transformation
$z=\frac{1+i\sigma}{1-i\sigma}$. That is, if $z$ is a unit eigenvalue of the
eigenproblem $Av=zBv$ then $\sigma$ is a real eigenvalue of the
eigenproblem $(A-B)v=\sigma (iA+iB)v$. 
The reason for this is that it is believed that finding real eigenvalues is an easier problem than
finding eigenvalues on the unit circle. Here, we apply a method doing
rational Krylov scans (c.f. \cite{Ruhe:1998:PRACTICAL}, \cite{Saad:1980:VARIATIONS}) along the real line implemented in the command
\texttt{sptarn} in the Matlab
control toolbox. However, we note that further exploitation of the
properties are necessary if we wish to treat large dimensional
DDEs.
An example of such a structure-exploitation is given in 
\cite{Jarlebring:2006:GAMM}.


\section{Examples}

\begin{example}{$n=35$, $m=3$}\label{ex:heatexample}
In order to determine how well the method scales with dimension we apply 
 the method in Section~\ref{sect:numerics}  to the discretization of a
 partial delay-differential equation. We
 consider the three delay partial
 differential equation,
\begin{equation}
 \left\{
\begin{split}
 u'_t(x,t)&=u''_{xx}(x,t)+\beta (1+\sin(3\pi x))u(x,t) \\
       &-\kappa_0\delta(x-x_0)u(0,t-h_1)\\
       &-\kappa_1\delta(x-x_1)u(x_1,t-h_2)\\
       &-\kappa_2\delta(x-x_2)u(1,t-h_3)\\
   u_x(0,t)&=0\\
   u_x(1,t)&=0, 
\end{split}
\right.
\label{example:pdde}
\end{equation}
where we pick $\kappa_0=\kappa_2=4$, $\kappa_1=10$, $x_0=1/3$,
 $x_1=1/2$, $x_2=3/4$ and $\beta=10$. The physical interpretation of 
 equation \eqref{example:pdde} is the heat equation on a rod with length one with
 heat production over the whole rod causing instability and three
  delayed  stabilizing pointwise feedbacks. This system is discretized with
 central difference in space with $n$ equally distributed intervals,
 yielding a system of the form 
 $\dot{x}(t)=A_0x(t)+A_1x(t-h_1)+A_2x(t-h_2)+A_3x(t-h_3)$, where 
 $x(t)\in \RR^n$, i.e., a three delay DDE. For this example we pick a
 step-length such that the discretized system is of dimension $n=35$.

 The critical surface closest to the origin of the discretized system is
 plotted in \fig{fig:heatexample}. The algorithm computes points
 on the critical surface, but for visualization purposes, we  connect the
 points to form a surface. 

 The application of the algorithm in Section \ref{sect:numerics} with the
 described exploitation for 390 different
 combinations of the parameters requires 62 minutes on a 
 computer running Linux and Matlab 7 on a $2.4$ GHz Intel Pentium 4 
 processor with 512 Mb RAM.


\begin{figure}[h]
  \centering
  \scalebox{0.6}{\includegraphics{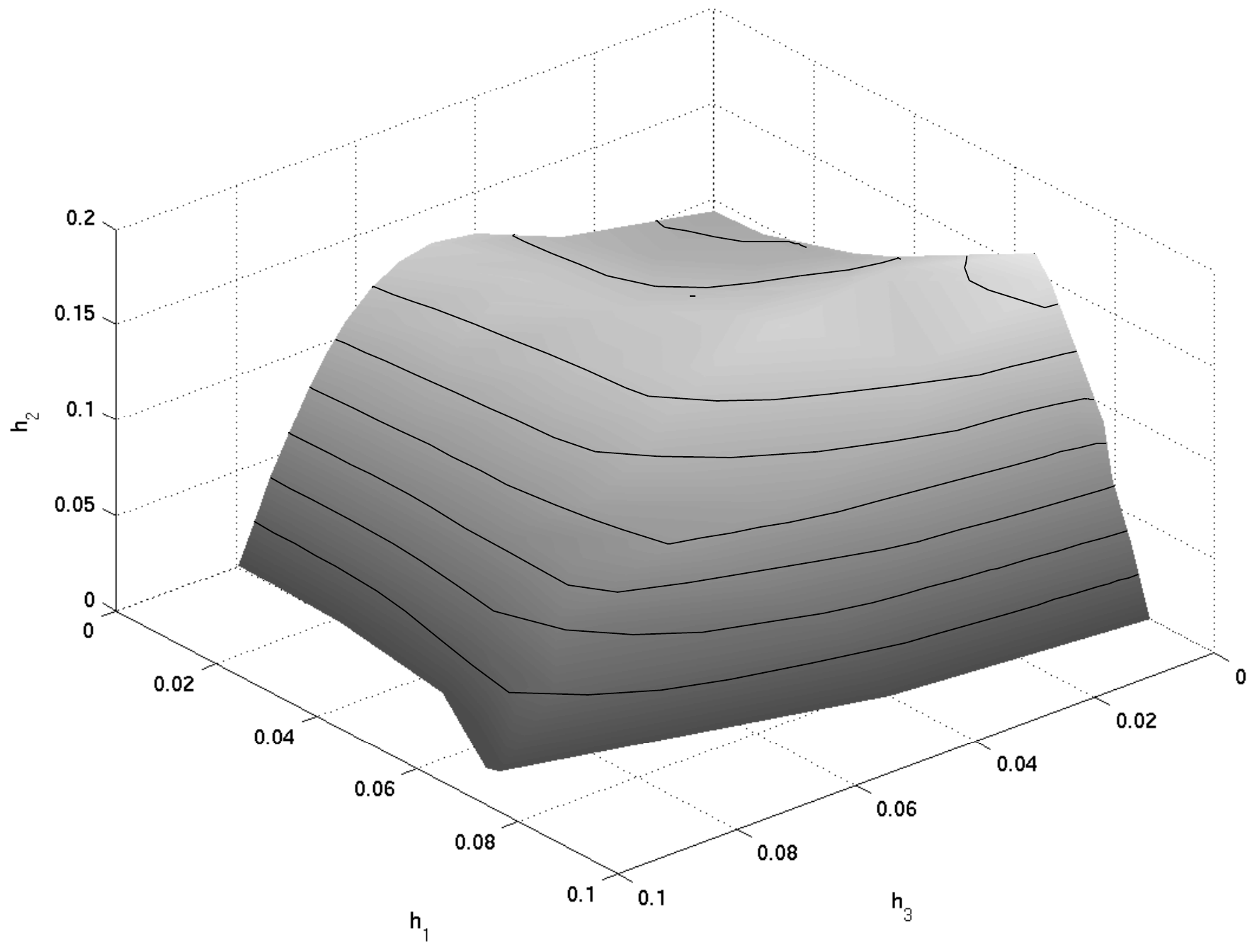}}
  \caption{Boundary of stability region for Example~\ref{ex:heatexample}}
 \label{fig:heatexample}
\end{figure}

\end{example}

\section{Conclusions}
We study the conditions on the delays for multiple delay time
delay systems such that the system has an imaginary eigenvalue. This is
done by introducing a condition using a function $\LL$, for which
the imaginary eigenvalue condition has a particularly easy form. For the
commensurate case the condition can be reduced to the
computation of the eigenvalues of a matrix similar to that of
Chen et al \cite{Chen:1995:DELAYMARGINS}. For the general (free) case we propose a
new method, which involves solving a quadratic eigenvalue problem. 

For the scalar case, the quadratic eigenvalue problem can be explicitly
solved and we can 
find a closed explicit expression for the condition on the delays. For
systems of large dimension we make initial remarks on what properties of
the problem an adapted
numerical scheme should exploit.

In the examples section we show the efficiency of the method by applying
it to previously solved examples as well a previously unsolved
three-delay example
of larger dimension.

\section*{Acknowledgments}
We would like to thank Heike Fassbender and Tobias Damm for very
valuable discussions  and comments.
\hfill





\bibliographystyle{elsart-num-sort}
\bibliography{fullbib}

\end{document}